\documentclass[11pt]{article}
\usepackage{graphicx}
\usepackage{amssymb}
\usepackage{epstopdf}
\usepackage{epsfig}
\usepackage{subfigure}
\usepackage{amsmath, amssymb, amsfonts}
\usepackage{a4, a4wide}
\usepackage{lscape}
\usepackage{color}
\usepackage{verbatim}
\usepackage{array}
\def\ie{{\rm i.e.,\/}\ }
\newcommand{\nco}{\newcommand}
\nco{\one}{\ensuremath{\,\,\mathrm{l}\!\!\!1}} 
\nco{\ZZ}{\mathbb{Z}}
\nco{\CC}{\mathbb{C}}
\nco{\red}{\color{red}}
\nco{\redend}{\normalcolor}
\nco{\0}{\underline{0}}
\newtheorem{theorem}{Theorem}

\newtheorem{definition}{Definition}
\begin{document}

\noindent
\begin{flushright}
{\small CERN-PH-TH-2010-060} \\
\end{flushright}
 
 \vspace{0.2cm} 
 
 %Title
\begin{center}
\begin{Large}
%\begin{bf}
    {Global dimensions for Lie groups at level $k$\\  \vspace{0.20cm} and their conformally exceptional quantum subgroups}
%\end{bf}
\end{Large}
\end{center}

%Authors

\begin{center}
\renewcommand{\thefootnote}{\arabic{footnote}}
\bf{R. Coquereaux} \footnotemark[1] ${ }^{ , }$\footnotemark[2]
\\ \renewcommand{\thefootnote}{\arabic{footnote}} % this has to be repeated
\end{center}

 \vspace{0.2cm} 
 
\begin{center}
{March 2010}
\end{center}

 \vspace{0.3cm} 

\abstract{We obtain formulae giving global dimensions for fusion categories defined by Lie groups $G$ at level $k$ and for the associated module-categories obtained via conformal embeddings. The results can be expressed in terms of Lie quantum superfactorials of type $G$. The later are related, for the type $A_r$,  to the quantum Barnes function.}

\vspace{0.3cm}


\vspace{0.2cm}

\noindent{\bf{Classification}}: 81R50; 81T40; 18D10; 33E99. 
% For 2000 Mathematics Subject Classification see http://www.ams.org/msc/

\addtocounter{footnote}{0}
\footnotetext[1]{ {\scriptsize{\it CERN, Geneva, Switzerland}. On leave from {\it CPT, Luminy, Marseille, France}.}}
\footnotetext[2]{\scriptsize{Talk given at: Coloquio de \'Algebras de Hopf, Grupos Cu\'anticos y Categor\'ias Tensoriales, C\'ordoba, Argentina, 2009}}

\vspace{0.4cm}

\unitlength = 1mm

%%%%%%%%%%%%%%%%%%%
\section{Introduction}

\subsection{Purpose}

To every complex Lie group $G$ and to every positive integer $k$ (the level), one associates a fusion category ${\mathcal A}_k={\mathcal A}_k(G)$, which is modular, and which is usually described in terms of integrable representations of affine Lie algebras, or in terms of a particular class of representations of quantum groups at roots of unity.  These fusion categories somehow generalize, at the quantum level, the theory of representations of Lie groups and Lie algebras. In classical group theory, groups may have non trivial subgroups and  the space of characters of a subgroup is a module over the ring of characters of the group.  In the present situation, we have a fusion ring, the Grothendieck ring of the monoidal category  ${\mathcal A}_k$, and one may consider quantum analogs of subgroups (or modules) by looking at module-categories ${\mathcal E}$ over which the given fusion category acts.  The purpose of the present article is mostly to present general formulae giving global dimensions for categories such as ${\mathcal A}_k$ and for a particular type of associated module-categories, that we call ``conformally exceptional quantum subgroups of Lie groups''.  In the process, we are led to define Lie quantum superfactorials  for all simple Lie groups. If $G=A_r$ their values can be obtained
from the quantum Barnes function when the argument is a  positive integer and the  $q$-parameter is an appropriate root of unity.
From the theoretical physics point of view, the results can be interpreted in the framework of boundary conformal field theories (WZW-models), or in the discussion of Chern-Simons topological field theories defined over $S^3$.

\subsection{Structure of the paper}

In the next subsection, we set the stage, summarizing  the necessary information.  In section 2 we consider fusion categories of type ${\mathcal A}_k(G)$ and obtain formulae for their global dimensions. The argument goes as follows: 1) the global dimension can be obtained from the square inverse of a particular element of the modular matrix $S$ implementing the $SL(2,\ZZ)$ action on this modular category, 2) this matrix element is  proportional to the quantum Weyl denominator, 3) The later is found, by inspection, to be equal to a product of quantum factorials over the exponents of $G$ (some care has to be taken in the non simply laced cases). 
In section 3 we consider associated module-categories, called ${\mathcal E}$.  
Thinking of $\mathcal{A}_k(G)/\mathcal{E}$ as a homogenous space, both discrete and quantum, we discuss a quantum analog of the Peter-Weyl theorem leading to an expression relating the global dimensions of $\mathcal{A}_k(G)$ and $\mathcal{E}$ to the dimension of  a particular subcategory  (called ambichiral) of $End_{{\mathcal A}_k(G)}{\mathcal E}$. 
When ${\mathcal E}$ is obtained from a conformal embedding, the ambichiral dimension is particularly simple to calculate. We therefore restrict our framework to those cases in order to obtain explicit formulae for global dimensions. Several remarks relating those results to geometry and physics are gathered in the last section.

\subsection{Background : General framework}
\label{background}
As discussed in  \cite{EtingofOstrik, FuchsSchweigert, KirilovOstrik, Ostrik}, action of the monoidal category ${\mathcal A}_k$ on the category ${\mathcal E}$ (assumed to be additive, semisimple and indecomposable but not necessarily monoidal) is obtained when we are given a monoidal functor from ${\mathcal A}_k$ to the monoidal category of endofunctors of ${\mathcal E}$.
Simple objects $a,b, \ldots$ of the module-category  ${\mathcal E}$ can also be thought  as right modules over a Frobenius  algebra  ${\mathcal F}$, which is a particular object in the monoidal category ${\mathcal A}_k$. 
From now on, ${\mathcal A}_k(G)$ is defined by the choice of a Lie group $G$ at level $k$  (see \cite{KazhdanLusztig}).  Using affine algebras or quantum groups is not necessary since all the tools that we need can be discussed within the framework of the theory of usual  finite dimensional Lie groups (see \ref{Liegroupslevelk}).
The reader only interested in the calculation of global dimensions of Lie groups at level $k$ may jump directly to the next section, but we gather below some standard material that will be used later.

The fusion ring of ${\mathcal A}_k$ comes with a special basis corresponding to simple objects $m,n,\ldots$, and its structure constants (non - negative integers) are encoded by the so - called fusion matrices $(N_n)^p_q$. 
 The  Grothendieck  group of ${\mathcal E}$ is a module over the Grothendieck ring of ${\mathcal A}_k$, and it is automatically a $\ZZ_+$ module: its structure constants are non negative integers often called annulus coefficients in string theory; they are encoded by ``annular matrices'' $F_n = (F_n)_{ab}$, where $a,b,\ldots$ refer to the simple objects of  ${\mathcal E}$ .
  The rigidity property of  ${\mathcal A}_k$ implies that the module ${\mathcal E}$ is rigid (or based \cite{Ostrik}). In other words:  $ (F_{\overline n})_{ab} = (F_ n)_{ba}$. 
The fusion graphs of  ${\mathcal A}_k$ are Cayley graphs describing multiplication by the generators (the corresponding matrices $N_f$ are their adjacency matrices).   The fusion graphs of  ${\mathcal E}$ are Cayley graphs describing the module action of the generators (the corresponding matrices $F_f$ are their adjacency matrices).  It is often useful to introduce the rectangular ``essential matrices'' $E_a = (E_a)_{nb} = (F_n)_{ab}$. 

Given a module-category ${\mathcal E}$ over a monoidal category  ${\mathcal A}_k$, we consider its endomorphism category  ${\mathcal O(E)} = End_{{\mathcal A}} {\mathcal E}$, which is also monoidal (call $x,y,\ldots$ its simple objects).
In practice one prefers to think in terms of rings and modules and use the same notations to denote them:  given a quantum module ${\mathcal E}$ over the fusion ring ${\mathcal A}_k$, one obtains an algebra ${\mathcal O(E)}$ ``of quantum symmetries'' (or double triangle algebra DTA \cite{Ocneanu:paths}) which  acts on ${\mathcal E}$ in a way compatible with the ${\mathcal A}_k$ action. 
The generators of ${\mathcal A}_k(G)$, \ie  those fundamental representations of $G$ that exist at the given level, also act on ${\mathcal O(E)}$, and in particular on its simple objects,  in two possible ways
\cite{Longo-Rehren, Evans-II}. The later is a bimodule over the former and one can associate to each of them two matrices with non negative entries describing left and right multiplication, and therefore two graphs, called left and right chiral graphs whose union is a non connected graph called the Ocneanu graph. The chiral graphs associated with one or another fundamental object of ${\mathcal A}_k$ can themselves be disconnected. 
More generally, the bimodule structure of  ${\mathcal O(E)}$ is described by ``toric matrices'' $W_{xy}$ defined by $m\, x \, n = \sum_y (W_{xy})_{m,n} \, y$. In particular, to 
every quantum module ${\mathcal E}$ is associated a matrix $Z= W_{00}$, where $0$ is the unit object of ${\mathcal O(E)}$;  it is such that $Z_{00}=1$ and it has the special property of being modular invariant, \ie it commutes with the action of the modular group $SL(2,\ZZ)$. This particular toric matrix (or its associated sesquilinear form) is often called in physics ``the modular invariant'', or  ``the partition function''.  If  ${\mathcal E} = {\mathcal A}_k $ it is just the identity matrix.  The other toric matrices  $W_{xy}$  can be physically interpreted in terms of  twisted partition functions (presence of defects in boundary conformal field theories, see \cite{PetkovaZuber:Oc}).  Two distinct quantum modules ${\mathcal E}$ sometimes lead to  the same partition function,  the classification issues are therefore different. 
Finally,  following the lines of  \cite{Ocneanu:paths}, see also \cite{Ostrik},  one can associate a quantum groupo\"\i d ${\mathcal B}$ to every module-category ${\mathcal E}$ over ${\mathcal A}_k(G)$. It is  a finite dimensional weak Hopf algebra which is simple and co-semisimple.

\section{Global dimensions for Lie groups at level $k$}

\subsection{Lie groups at level $k$} \label{Liegroupslevelk}
\subsubsection*{Simple objects of ${\mathcal A}_k(G)$ }

Although usually defined in terms of affine Lie algebra or quantum groups concepts, they can be simply obtained as follows.

{\definition{} Let $n$ denote an irreducible representation of the simple Lie group $G$, with highest weight $\lambda_n$. We define its level as the integer $\langle \lambda_n,  \theta  \rangle$ where  $\theta$ is the highest root of $Lie (G)$ and $\langle .. , .. \rangle $ is the fundamental quadratic form. We shall say that  $n$ is integrable at level $k$ (or that ``it exists at level $k$'') if and only if its level is smaller than $k$ or equal to $k$.} 

\smallskip
 The condition of existence at level $k$ selects a finite set, with cardinal $r_A$, of irreducible representations $m,n\ldots$ of $G$ that are parametrized by $r$-tuples in the basis of fundamental weights, $r$ being the rank of $G$.   This set can be conveniently ordered  by increasing values of the level; at a given level  the $r$-tuples are sorted in a standard way, starting from the end (see examples later).
  In the case $G=SU(N)$ all fundamental representations already appear at level $1$. This is usually not so for other choices of $G$. 
We call ``basic representations'' those fundamental representations of $G$ that have smallest classical dimension. 
Quantum dimension (or ``categorical dimension) $qdim(n)$ of a simple object $n$ can be calculated from the quantum version of the Weyl formula, together with the choice of a root of unity $q=exp(i \pi/\kappa)$, for an altitude $\kappa=k+g$, where $g$ is the dual Coxeter number of $G$. The order of ${\mathcal A}_k$, or global dimension, is $\vert {\mathcal A}_k \vert = \sum_n (qdim(n))^2$. It is sometimes called the mass of $G$ at level $k$ and denoted $\mathrm{mass}(G,k)$.

\subsubsection*{Action of $SL(2,\ZZ)$ and the  Kac-Peterson formulae}
 \label{Kac-Peterson}
 
Action of the group $SL(2,\ZZ)$ on the vector space spanned by the simple objects $m,n\ldots$ of the modular category ${\mathcal A}_k(G)$ is described by matrices $S$ and $T$ representing the two generators $\tau \mapsto -1/\tau$ and $\tau \mapsto \tau +1$ of the modular group. This representation (actually a representation of the double cover), known by Hurwitz  long ago \cite{Hurwitz} for $G=SU(2)$, is given, in general,  by the following formulae \cite{Kac-Peterson}:   
$$
S_{mn}= \frac{i^{\Sigma_{+}} \sqrt{\Delta} }{ (g+k)^{r/2} }\left(\sum _{w} \, \epsilon_{w}  \, e^{-\dfrac{2 i \pi  \langle w(m+\varrho ),n+\varrho \rangle }{g+k}}\right) 
,\;
T_{mn}=e^{2 i \pi \left[\dfrac{\langle m+\varrho ,m+\varrho \rangle }{2 (g+k)}-\dfrac{\langle \varrho ,\varrho \rangle }{2 g}\right]} \; \delta_{mn}
$$
where $g$ is the dual Coxeter number,  $k$ is the level, $w$ runs over the Weyl group of $G$, $\epsilon_{w}$ is its signature, $r$ is the rank of $G$, $\varrho$ is the Weyl vector, $\Sigma_{+}$ is the number of positive roots (also equal to the sum of exponents, \ie to $r \, \gamma/2$ where $\gamma$ is the Coxeter number),  and $\Delta$ is the
determinant of the fundamental quadratic form (also the inverse of the ``long index'' \ie the index, in the weight lattice, of the sublattice of long roots). 
With these definitions one has $(S\,T)^3 = S^2 = \mathcal{C}$, the ``charge matrix'' satisfying $\mathcal{C}^2=\one$.
$S$ is unitary and symmetric.
Remember that $T$ is related as follows to other group theoretical quantities:
 The  eigenvalue of the quadratic Casimir for a representation $\lambda$ is $C(\lambda)=\langle \lambda ,\lambda +2 \varrho \rangle$, the classical Dynkin index is $dim(\lambda)\, C(\lambda)/2 d$, where $d=dim(G)$, and, at level $k$, the conformal weight of $\lambda$ is defined as\footnote{The coefficient $2$ stands for $\langle \lambda ,\lambda \rangle$, the usual convention for long roots.} $h(\lambda)=C(\lambda)/2(k+g)$. Some authors, like \cite{BakalovKirillov}, prefer to use a matrix $t$ that differs from the one introduced previously by a modular phase: $t = T \, exp(2 i \pi c/24)$ where $c$ is the central charge $c={d \, k}/{(g+k)}$.

 \subsubsection*{Intertwining operators (morphisms) and the Verlinde formula}
 
 Simple objects of the category ${\mathcal A}_k(G)$ have been discussed in the previous subsection. 
 There is a morphism from $m \otimes n$ to $p$ whenever $ (N_m)_{np} \neq 0$. 
 One way to obtain the morphisms, in practice, is to use the Verlinde formula. We recall it below for the convenience of the reader. It expresses the fusion matrices in terms of the  modular generator $S$ of $SL(2,\ZZ)$  :
    \label{Verlinde}
    $$ (N_m)_{np} = \sum_q \, \frac{S_{mq} S_{nq} S_{pq}^\star}{S_{0q}}$$
where $m,n,p$ are simple objects and $0$ refers to the trivial object.

\subsection{Weyl formula, $S$ matrix and quantum dimensions}

\subsubsection*{The quantum Weyl formula and the quantum Weyl denominator}

If a representation $n$ exists at level $k$, its quantum dimension at level $k$ is given by the quantum version of the Weyl formula, given below, with $q=exp(i \pi/(g+k)$. The classical dimension is obtained by taking $q=1$. This formula can be obtained from the theory of quantum groups or from the theory of affine Lie algebras, but in our framework we take this formula as a definition.
 Supposing that $n$ is irreducible, we  use the same notation for the representation and for its highest weight.
\[qdim (n)= \prod_{\alpha > 0} \frac{\langle n + \rho, \alpha\rangle_q}{\langle\rho, \alpha\rangle_q}\] where $\rho$ is the Weyl vector and $\alpha$ belongs to the set of positive roots.
There are two standard definitions for q-numbers:  $n_q = ~ [n]_q$, that we used in the above formula, and $[[n]]_q$:
$$ [n]_q = \frac{q^n - q^{-n}}{q-q^{-1}} \quad \text{and}  \quad [[n]]_q =  \frac{1 - q^{n}}{1-q}$$
The quantum Weyl denominator defined below would be different if we had made the other choice. 
Using $q$-numbers $[n]_q$, and with the same notations as before we have :
{\definition{} The quantum Weyl denominator of a Lie group $G$ at level $k$ (a non negative integer), is the number ${\mathcal D}_W=\prod_{\alpha > 0} \,   \langle\rho, \alpha\rangle_q$. }

\subsubsection*{$SU(2)$-fusion and the ribbon diagram of $G$}

It is convenient to display the scalar products between a chosen weight  and all the roots of $G$ in the Ocneanu ribbon diagram (seen as a generalized root set  \cite{Ocneanu:MSRI}) 
associated with the Dynkin diagram of $G$.
It is a periodic strip of size $r \times \gamma$ (actually a $Z_2$ quotient of a rectangle $(2r) \times \gamma$), with vertical period equal to $2 \gamma$.

If $G$ is simply laced ($g=\gamma$), and since-module categories of type $SU(2)$ are classified by $ADE$ Dynkin diagrams \cite{CIZ, KirilovOstrik},
 to every vertex $a$ of its  Dynkin diagram, with $r$ vertices and Coxeter number $\gamma$,  one can associate both a fundamental weight $\omega_a$ of a Lie group $G$ and a simple object $a$ of an  $SU(2)$ module-category at level $\nu = \gamma -2$ (this $2$ stands for the Coxeter number of $A_1 \sim SU(2)$).
The Dynkin diagram, that we call also $G$, like its adjacency matrix, is the fundamental fusion graph of the later. The adjacency matrix determines both the fundamental quadratic form  (inverse of the Cartan matrix $2 \, \one   - G$) of the root system, and the family of all annular matrices $F_n=(F_n)_{a,b}$ labelled by simple objects $n$ of ${\mathcal A}_\nu(SU(2)$, identified with integers, thank's to the usual Tchebychef recurrence formula $F_n = F_{n-1} G - F_{n-2}$, $F_0= \one$, $F_1=G$. When $\gamma$ is finite this defines a periodic family of $r \times r$ matrices, with period $2\gamma$. 
It is convenient to define essential (rectangular) matrices $E_a = (E_a)_{nb} = (F_n)_{ab}$. The key observation \cite{Ocneanu:MSRI} is that because of the ubiquitous role of the adjacency matrix $G$,  matrix elements $(E_a)_{nb}$ can be interpreted either as scalar products between the fundamental weight $\omega_a$ and all the $r \, \gamma$ roots, or as fusion coefficients (also as the number of essential paths of length\footnote{In the case of $SU(2)$ irreducible representations are labelled by a single integer.} $n$ from $a$ to $b$ on the graph $G$). More generally, the same ribbon can be used to display and calculate scalar products between arbitrary weights and all the roots.

If $G$ is not simply laced, one can still draw the ribbon diagram, but one has to be careful since the ribbon of roots and the ribbon of coroots are not equal (the columns associated with short roots differ from the same columns for coroots by constant scaling coefficients).
To illustrate this we display below the scalar products between an arbitrary weight $\omega = \sum_a \lambda_a \, \omega_a$ of $A_5$, of $D_4$ and of $G_2$ and all the positive roots of the chosen Lie groups (so only half of the Ocneanu ribbon is drawn). We also give the Weyl vector $\rho$. For all three cases, $\gamma=6$, but $g=6$ for $A_5$ and  $D_4$ whereas $g = 4$ for $G_2$.  The obtained values for scalar products are immediately obtained from fusion matrices (for $A_5$) or annular matrices (for $D_4$) describing respectively the fusion category ${\mathcal A}_4(SU(2))$ and the $SU(2)$ module-category with fusion graph $D_4$. In the case of $G_2$ one has to rescale the last column. The three cases are obtained from the study of $SU(2)$ at level $\nu=4$ 
since $\nu=\gamma-2$.

Case $A_5$:

\smallskip
$
{\scriptsize
\omega = 
\begin{array}{lllll}
 \lambda_1 &  & \lambda_3 &  & \lambda_5 \\
  & \lambda_1+\lambda_2+\lambda_3 &  & \lambda_3+\lambda_4+\lambda_5 &  \\
 \lambda_2+\lambda_3 &  & \lambda_1+\lambda_2+\lambda_3+\lambda_4+\lambda_5 &  & \lambda_3+\lambda_4 \\
  & \lambda_2+\lambda_3+\lambda_4+\lambda_5 &  & \lambda_1+\lambda_2+\lambda_3+\lambda_4 &  \\
 \lambda_4+\lambda_5 &  & \lambda_2+\lambda_3+\lambda_4 &  & \lambda_1+\lambda_2 \\
  & \lambda_4 &  & \lambda_2 & 
\end{array}
}
$

\smallskip
$
\begin{array}{cc}

{\scriptsize
\rho=
\begin{array}{lllll}
 1 &  & 1 &  & 1 \\
  & 3 &  & 3 &  \\
 2 &  & 5 &  & 2 \\
  & 4 &  & 4 &  \\
 2 &  & 3 &  & 2 \\
  & 1 &  & 1 & 
\end{array}
}

&

\unitlength 0.5mm
\begin{picture}(5,5)
%%%%%%
\put(20,20){\line(1,-1){10}}
\put(40,20){\line(-1,-1){10}}
\put(40,20){\line(1,-1){10}}
\put(60,20){\line(-1,-1){10}}
%%%%%%
\put(50,10){\circle*{2}}
\put(30,10){\circle*{2}}

\put(20,20){\circle*{2}}
\put(40,20){\circle*{2}}
\put(60,20){\circle*{2}}
%%%%%
\end{picture}

\end{array}
$

\smallskip
Case $D_4$:

\hskip 8.5cm
\unitlength 0.5mm
\begin{picture}(5,5)
%%%%%%
\put(20,20){\line(1,-1){10}}
\put(40,20){\line(-1,-1){10}}
\put(30,10){\line(3,1){30}}
%%%%%%
\put(30,10){\circle*{2}}

\put(20,20){\circle*{2}}
\put(40,20){\circle*{2}}
\put(60,20){\circle*{2}}
%%%%%
\end{picture}

\smallskip
$
{\scriptsize
\omega = 
\begin{array}{llll}
 \lambda_1 &  & \lambda_3 & \lambda_4 \\
  & \lambda_1+\lambda_2+\lambda_3+\lambda_4 &  &  \\
 \lambda_2+\lambda_3+\lambda_4 &  & \lambda_1+\lambda_2+\lambda_4 & \lambda_1+\lambda_2+\lambda_3 \\
  & \lambda_1+2 \lambda_2+\lambda_3+\lambda_4 &  &  \\
 \lambda_1+\lambda_2 &  & \lambda_2+\lambda_3 & \lambda_2+\lambda_4 \\
  & \lambda_2 &  & 
\end{array}
},
\qquad 
{\scriptsize
\rho=
\begin{array}{llll}
 1 &  & 1 & 1 \\
  & 4 &  &  \\
 3 &  & 3 & 3 \\
  & 5 &  &  \\
 2 &  & 2 & 2 \\
  & 1 &  & 
\end{array}
}
$

\smallskip
Case $G_2$:

\smallskip
$
\begin{array}{ccc}

{\scriptsize
\omega = 
\begin{array}{ll}
 \lambda_1 &  \\
  & \lambda_1+\frac{\lambda_2}{3} \\
 2 \lambda_1+\lambda_2 &  \\
  & \lambda_1+\frac{2 \lambda_2}{3} \\
 \lambda_1+\lambda_2 &  \\
  & \frac{\lambda_2}{3}
\end{array}
},
& 
\qquad
{\scriptsize
\rho=
\begin{array}{ll}
 1 &  \\
  & \frac{4}{3} \\
 3 &  \\
  & \frac{5}{3} \\
 2 &  \\
  & \frac{1}{3}
\end{array}
}
&
\unitlength 0.8mm
\begin{picture}(5,5)
%%%%%%
\put(20,10){\line(1,-1){10}}
\qbezier(20, 10)(40,10)(30, 0)
%%%%%%
\put(12,2){$\nwarrow^3$}
\put(35,7){$\searrow^1$}
%%%%%%
\put(30,0){\circle*{2}}

\put(20,10){\circle*{2}}
%%%%%
\end{picture}
\end{array}
$

\bigskip
\noindent
In the case of $G_2$ for instance, the quantum dimension of a representation of highest weight $\{\lambda_1,\lambda_2\}$, in the basis of fundamental weights, is read from the last ribbon diagrams as 
$$
\frac{\left(\lambda _1+1\right)_q \left(\frac{\lambda _2}{3}+\frac{1}{3}\right)_q \left(\lambda _1+\frac{\lambda
   _2}{3}+\frac{4}{3}\right)_q \left(\lambda _1+\frac{2 \lambda _2}{3}+\frac{5}{3}\right)_q \left(\lambda
   _1+\lambda _2+2\right)_q \left(2 \lambda _1+\lambda _2+3\right)_q}{\left(\frac{1}{3}\right)_q 1_q
   \left(\frac{4}{3}\right)_q \left(\frac{5}{3}\right)_q 2_q 3_q}
$$ Notice that the Weyl denominator of $G$ is directly  obtained by multiplying all the entries of the ribbon diagram associated with the weight $\rho$.

  \subsubsection*{From the matrix element $S_{00}$ to the quantum Weyl denominator}
  
  Let $\alpha$ be a positive root, $\rho$ the Weyl vector, and $q=exp(i \pi/(g+k)$. The q-scalar product $\langle \rho, \alpha \rangle_q =(q^{\langle \rho, \alpha \rangle}-q^{-\langle \rho, \alpha \rangle})/(q-q^{-1})$, where $\langle . , . \rangle$ is the fundamental quadratic form, can be written
  $$\langle \rho, \alpha \rangle_q = (exp(\frac{\alpha}{2}) - exp(\frac{-\alpha}{2})) \, [\frac{2 i \pi}{k+g}] \, / \, (2 i \, \sin\frac{\pi}{g+k})$$
  where $e^\lambda$ is a formal exponential defined by $e^\lambda[\mu] = e^{\langle \lambda, \mu \rangle}$, for  $\lambda, \mu$,  arbitrary elements of the space of weights.
  A standard manipulation (see for instance \cite{YellowBook}) allows one to write 
  $$\prod_{\alpha >0} (e^{\alpha/2}-e^{-\alpha/2}) = \sum_{w \in W} \, \epsilon_w \, e^{w \rho}$$
  where $W$ is the Weyl group.
  Taking $m=n=0$, the trivial weight, in the Kac-Peterson formula for $S_{mn}$, and using the fact that the number of positive roots is $r \gamma/2$, we obtain the following (notice that the $i$ factors cancel).
  
 {\theorem {}
The element $S_{00}$ of the modular matrix $S$ and the quantum Weyl denominator $ {\mathcal D}_W = \prod_{\alpha >0} \,  \langle \rho, \alpha \rangle_q  $ are related as follows: 
$$
S_{00} = \frac {2^{\frac{r \gamma}{2}} \, \sqrt{\Delta}}
{{(g+k)}^{r/2}} \;  {\left( \sin \frac{\pi}{g+k} \right) }^{\frac{r \gamma}{2}} \; {\mathcal D}_W$$
}

{\sl Remark:} One should remember that there are two standard definitions for q-numbers,  $n_q = ~ [n]_q$, that we use in this paper, and $[[n]]_q$. The quantum Weyl denominator, together with the above relation,  look different if we make the other choice. The above relation between $S_{00}$ and ${\mathcal D}_W$, although written differently (for instance in \cite{YellowBook}), is, of course, hardly new.

\subsubsection*{Quantum dimensions and the $S$ matrix}

From the quantum Weyl formula, and performing standard manipulations, like those done above, one obtains the well known relation:
$$ qdim(n) = S_{n0}/S_{00}$$
Both the numerator and the denominator are positive real numbers.
Unitarity of the $S$ matrix implies $\sum_n \, \vert S_{n0} \vert^2 =1$. In particular, the global dimension $\vert {\mathcal A}_k(G) \vert$, defined as the sum of squares of the q-dimensions of its simple objects, is obtained as:

$$\sum_n \, qdim(n)^2 = \frac{1}{S_{00}^2}$$

\subsubsection*{Classical character polynomials and the $S$ matrix}

The following result is quite handy but does not seem to be so well known. We mention it without proof.
Call $\chi(m)=\chi(m; t_1,t_2,\ldots t_r)$ the classical character polynomial of the Lie group $G$, associated with an irreducible representation $m$ defined by its highest weight also denoted $m$. It encodes the weight system of $m$ :  each weight $p$ occurring with multiplicity $a$ in this weight system gives a term $a\, t_1^{p_1} t_2^{p_2} \ldots t_r^{p_r}$ in the polynomial. Here $(p_1, p_2, \ldots p_r)$ are the components of $p$ on the basis of fundamental weights.
We know how to evaluate a monomial of the type $ t_1^{p_1} \, t_2^{p_2} \ldots t_r^{p_r}$ (powers can be positive or negative) on a weight $n$ :
This value, by definition, is $exp[2 i \pi/(g+k) \; \langle p_1 \omega_1 + p_2 \omega_2 + \ldots + p_r \omega_r, n  \rangle]$ where $\omega_i$ are the fundamental weights. This evaluation is extended by linearity so that one can evaluate the character polynomial $\chi(m)$ on arbitrary weights $n$. The obtained value is denoted $\chi(m)[n]$. 
Assuming that $m$ and $n$ are two irreducible representations of $G$ existing at level $k$, one finds the following relation between the matrix elements of  $S$  and the character polynomial (remember that the Weyl group permutes the weights of a weight system) : 
$$S_{mn} / S_{00} =  qdim[n] \; \chi(m)[n + \rho]$$
It gives a way to calculate the matrix $S_{mn}/S_{00}$ from known classical character polynomials.
The above relation looks highly asymmetrical but implies  $qdim[n] \, \chi(m)[n+ \rho] = qdim[m] \, \chi(n)[m+ \rho]$ since $S$ is symmetric.
In particular $S_{m0} / S_{00} =  \chi(m)[ \rho]$. Now, since $\rho = \omega_1 + \omega_2 + \ldots + \omega_r$, where $\omega_j$ are fundamental weights, $\rho^j=\langle\omega_j, \rho\rangle$ are the components of the Weyl vector on the base of simple coroots and $exp(\frac{2 i \pi}{g+k} \, p_1 \, \langle \omega_j, \rho \rangle) = q^{2 p_1 \rho^j}$.
Assuming that $m$ exists at level $k$,  we have also  $$qdim(m) = \chi(m; q^{2 \rho^1},q^{2 \rho^2} ,\ldots q^{2 \rho^r})$$ where $\chi(m)=\chi(m; t_1,t_2,\ldots t_r)$ is the classical character polynomial of $m$ and  $(\rho^j)$ are the components of the Weyl vector on the base of simple coroots.

\subsubsection*{Quantum dimensions and the Perron-Frobenius vector}

In some cases the fusion graphs of the studied example are somehow given. In such a situation, it is usually simpler to calculate the quantum dimensions of simple objects from the Perron-Frobenius vector associated with adjacency matrices (fusion matrices) of these graphs, since it coincides, up to normalization, with the first line of the modular matrix $S$, thought as an analog of a character table for finite groups.  It is actually enough to consider the fusion graph relative to a basic representation.

\subsection{Global dimensions and Lie superfactorials}

\subsubsection*{Lie superfactorial and quantum Lie superfactorial}

{\definition{}Let $G$ be a simply laced complex simple Lie group and q be a complex number. Call ${\mathcal E}$ the multiset of exponents of $G$. We define the 
quantum Lie superfactorial of type $G$ as follows:
$$ \mathrm{sf}_{G}[q]= \prod_{s \in {\mathcal E}} \,  [s]!_q $$
}

\noindent
Here $[s]!_q$ denotes the quantum factorial of $s$, \ie 
$$[s] !_q = \prod_{n=1}^{n=s} [n]_q \quad \text{with} \quad  [n]_q = \frac{q^n - q^{-n}}{q-q^{-1}}$$

  \noindent
{\sl Expressions for the quantum Lie superfactorial}:
\begin{itemize}
\item
  If  $G=A_{r}\sim SU(n+1)$, then ${\mathcal E} = \{1,2,\ldots,r\}$, and $ \mathrm {sf}_q(r) \doteq \mathrm{sf}_{A_r}[q]= \prod_{s=1}^{s=r} \,  [s]!_q $.
\newline  \noindent
 Notice that  $  \mathrm{sf}(r) \doteq \lim_{q \to 1}  \mathrm{sf}_{A_r}[q] =  \prod_{s=1}^{s=r} \, s! $  
\newline  \noindent
The classical limit of the quantum Lie superfactorial of type $A_r$ therefore coincides with the superfactorial function defined in \cite{SloanePlouffe}.

\item  If  $G=D_{r}\sim SO(2r)$, set $g=2r-2$, then ${\mathcal E} = \{1,3,5,\ldots g-3, g-1; g/2\}$ (when $r$ is even, $g/2$ appears twice), and $ \mathrm{sf}_{D_r}[q]= [g/2]!_q \, \prod_{s=1,3,5,\ldots}^{s=g-1} \,  [s]!_q $.
  
\item  If $G = E_6$ then $ \mathrm{sf}_{E_6}[q]=[1]!_q  \, [4]!_q  \, [5]!_q \, [ 7]!_q  \,  [8]!_q  \,  [11]!_q $
\newline  \noindent
 Its classical limit is the integer $\mathrm{sf}_{E_6}={{2^{25}} {3^{10}}{5^5} {7^3} {11^1}}$.
  
\item  If $G = E_7$  then $ \mathrm{sf}_{E_7}[q]=[1]!_q  \, [5]!_q   \,  [7]!_q  \,  [9]!_q  \,  [11]!_q \, [13]!_q \, [17]!_q $
\newline  \noindent
  Its classical limit is the integer $\mathrm{sf}_{E_7}=2^{47} 3^{22} 5^{10} 7^{6} 11^{3} 13^{2} 17^{1}$.

\item   If $G = E_8$  then   $\mathrm{sf}_{E_8}[q]=[1]!_q  \, [7]!_q \,  [11]!_q \,  [13]!_q \, [17]!_q \,  [19]!_q \,  [23]!_q \,  [29]!_q$
\newline  \noindent   
    Its classical limit is the integer   $\mathrm{sf}_{E_8}= 2^{97} 3^{47} 5^{21} 7^{14} 11^{8} 13^{6} 17^{4} 19^{3} 23^{2} 29^{1}$.
   
  \end{itemize}

  {\definition{}Let $G$ be a non-simply laced complex simple Lie group and q be a complex number. Call ${\mathcal E}$ the multiset of exponents of $G$. We define the 
quantum Lie superfactorial of type $G$ as follows:
\begin{itemize}
\item If $G=B_{r} \sim SO(2r+1)$, then ${\mathcal E} = \{1,3,5,\ldots ,2r-1\} $
$$ \mathrm{sf}_{B_r}[q]= \prod_{s \in {\mathcal E}} \,  \widetilde{[s]!_q}  \quad \text{with} \quad  \widetilde{[s]!_q}  = [\frac{s}{2}]_q \,  [s-1]_q  \,  [s-2]_q \ldots [3]_q  \, [2]_q \,  [1]_q,  \quad \text{and} \quad   \widetilde{[1]!_q} =[1/2]_q$$ 
Its classical limit is the rational $ \mathrm{sf}_{B_r} = \frac{1}{2^r} \prod_{s, odd=1}^{2r-1} s!$ 
\item If $G=C_{r} \sim Sp(2r)$, then  ${\mathcal E} = \{1,3,5,\ldots ,2r-1\} $
 $$ \mathrm{sf}_{C_r}[q]= \prod_{s \in {\mathcal E}} \, \widetilde{[s]!_q}  \quad \text{with} \quad  \widetilde{[s]!_q}  = \left[\frac{s}{2}\right]_q \, \left[\frac{s-1}{2}\right]_q \ldots
 \, \left[\frac{s-\frac{s-3}{2}}{2}\right]_q \,  \left[{s-\frac{s-1}{2}}\right]_q \, \left[\frac{s-\frac{s+1}{2}}{2}\right]_q  \ldots \left[\frac{2}{2}\right]_q \,   \left[\frac{1}{2}\right]_q$$
 $ \text{and} \quad   \widetilde{\left[1\right]!_q} =1$.
Its classical limit is the rational $ \mathrm{sf}_{C_r} = \dfrac{1}{2^{r(r-1)}} \prod_{s, odd=1}^{2r-1} s!$ 
\item If $G=F_4$
$$ \mathrm{sf}_{F_4}[q]=\left[\frac{1}{2}\right]_q^2 1_q^3 \left[\frac{3}{2}\right]_q 2_q^3
   \left[\frac{5}{2}\right]_q^2 3_q^3 \left[\frac{7}{2}\right]_q 4_q^2
   \left[\frac{9}{2}\right]_q 5_q^2 \left[\frac{11}{2}\right]_q 6_q 7_q
   8_q $$
    Its classical limit is the integer $ \mathrm{sf}_{F_4} = (1/2^{12}) \, 1! \, 5! \, 7! \, 11! = 2^{15} 3^7  5^ 4 7^2 11^1$
\item If $G=G_2$
$$ \mathrm{sf}_{G_2}[q]=  \left[\frac{5}{3}\right]_q \left[\frac{4}{3}\right]_q 3_q \,   2_q \,  1_q \,   \left[\frac{1}{3}\right]_q$$   
 Its classical limit is the rational $ \mathrm{sf}_{G_2} =  \dfrac{1! \,  5! }{ 3^3}=   40/9$.
\end{itemize}
}

 \noindent
{\sl Remarks}:
\begin{itemize}
\item
{\sl Warning} : There are two non-equivalent standard definitions for $q$-numbers, and, correspondingly, two non-equivalent definitions for the quantum factorials and the quantum Lie superfactorials of type $G$. The other possibility\footnote{$[[s]]!_q$ and $(a;q)_s$ coincides with the functions $QFactorial[s,q]$ and $QPochhammer[a, q, s]$ of {\sl Mathematica}}  (using a notation with double bracket) is indeed to set:
$$ [[s]]!_q =  \prod_{n=1}^{n=s} [[n]]_q  \quad  \text{with}  \quad [[n]]_q =  \frac{1 - q^{n}}{1-q}$$
We shall use mostly the first, namely $[s]_q$ rather than $[[s]]_q$ because it is real when $q$ is a root of unity. 
The two factorials related as follows  can be used for non integer values of their argument by using the  q-Pochhammer symbol $(a;q)_n= (a;q)_\infty / (a q^n;q)_\infty$, with $(a;q)_\infty = \prod_{k=0}^\infty \, (1-a q^k)$.
$$[s]! _q = q^{-\frac {s (s - 1)} {2}}\, [[s]]! _ {q^2} = q^{-\frac {s (s - 1)} {2}}\,  (1 - q^2)^{-s} \,  (q^2, q^2)_s $$
For quantum superfactorials of type $G$, we shall use mostly $\mathrm {sf}_G[q]$, as defined above. However, it may be convenient, sometimes, to use $[[n]]_q$  rather than $[n]_q$. For this reason we introduce the definition and notation $\mathrm {Sf}_G[q]$, with a capital S, by just replacing $[s]!_q$ by $[[s]]!_q$ in the above definitions. For instance, $ \mathrm {Sf}_q(r) \doteq \mathrm{Sf}_{A_r}[q]= \prod_{s=1}^{s=r} \, [[s]]!_q$. Notice that $\mathrm {Sf}_{q^2}(r) = q^{((r+1)r (r-1))/6}\,  \mathrm {sf}_{q}(r)$.
\item The definition of the quantum superfactorial, for the non-simply laced cases, has been chosen in such a way that it agrees with the value of the  quantum Weyl denominator.
\item Warning: At the classical level, the correcting factors (rational numbers) appearing, for non-ADE Lie groups, in the numerator and in the denominator of the Weyl formula, cancel out. One should be tempted of defining  classical Lie superfactorials for non-ADE Lie groups as the product of factorials of their exponents, like for the ADE cases, however,  such functions would do not coincide with the $q \mapsto 1$ limit of their quantum counterparts.
\end{itemize}

\bigskip

%%%%%%%%%

{\theorem{} The quantum Weyl denominator of $G$ is equal to the quantum superfactorial of type $G$. We have :  ${\mathcal D}_W = \mathrm{sf}_G[q]$.}

\smallskip
In the case $G=A_n$ the fact that the denominator is a product of successive factorials is a familiar result, often used when manipulating Young tableaux  together with the Schur dimension formula to calculate dimensions of irreducible representations. Writing the result in terms of exponents is immediate since they are consecutive integers $\{1,2,\ldots n\}$. For the other $ADE$ Lie groups, the theorem is obtained by inspection, namely from the explicit determination of the Weyl denominator in all cases (we found convenient to do these calculations by using the ribbon diagram of $G$ introduced in a previous section).  For the non-simply laced cases, the theorem results from our definition of quantum superfactorials involving non trivial pre-factors, since they were precisely chosen to match the calculated Weyl denominator.

\subsubsection*{Global dimensions and quantum Lie superfactorials}

Putting all previous results together, we obtain

{\theorem{} Let $G$ a simple Lie group with rank $r$, Coxeter number $\gamma$, dual Coxeter number $g$, and $k$ a non-negative integer.  The global dimension of the fusion category ${\mathcal A}_k(G)$ is given explicitly by 
$$ \vert {\mathcal A}_k(G)\vert \doteq \sum_n \, qdim(n)^2 = 
\frac{(g+k)^r}{2^{r \gamma} \, \Delta \,  {\left( \sin  \dfrac{\pi}{g+k}\right)}^{r \gamma} \, \left( \mathrm{sf}_G [q]\right)^2}
$$ 
where $sf_G [q]$ is the quantum Lie superfactorial of type $G$, $\Delta$ is the determinant of the fundamental quadratic form, and $q = \exp{i \pi/(g+k)}$.}

\subsubsection*{Classical limits} \label{classical-limit}

When $q \mapsto 1$ (or $k \mapsto \infty$), the Lie quantum superfactorials go to their classical limits, the Lie classical superfactorials (remember that for $ADE$ cases, the later are just equal to the product of the corresponding exponents, but they involve non trivial pre-factors for non-$ADE$ cases).  Remember also that $r+r\gamma$ is equal to the dimension of the Lie group $G$.  Therefore, 

$$ \vert {\mathcal A}_k(G)\vert \,  \underset{{{k\rightarrow \infty}} }{\sim}\, \frac{1}{2^{r\gamma} \, \Delta \, \pi^{r\gamma} \, \left({\mathrm{sf}_G}\right)^2}\times k^{dim_G}$$

It is pleasant to notice that the coefficient in front behaves like the square of the volume of $G$. Indeed, as shown by \cite{macdonald}, the later is given, up to some normalization factor, by a product of volumes $2 \pi^{n + 1}/n! $ of the spheres of  dimensions $2 n + 1$, $n$ running over exponents, that enter the De Rham cohomology of $G$.

\subsubsection*{On the classical and quantum Barnes G-function}

The classical Barnes G-function is an entire function that obeys the  functional equation $G(z+~1) = \Gamma(z) \, G(z)$.  For $n=1,2,3,\ldots$ it is related to the superfactorial: $G(n+2)=\mathrm {sf}(n)$. The functions $z$, $\Gamma(z)$ and $G(z)$ are actually the first three members of a hierarchy of functions known as Barnes multigamma functions $G^{(d)} (z)$ satisfying $G^{(d)} (1)=1$ and $G^{(d)} (z+1) = G^{(d-1)} (z) \, G^{(d)} (z)$.  
The first three members are $G^{(0)} (z)=z$, $G^{(1)} (z)=\Gamma(z)$ and $G^{(2)} (z)=G(z)$. Quantum analogs of these functions have been defined (see  \cite{Koshkin},  \cite{Nishizawa} for details), they also obey $G_q^{(d)} (1)=1$ and $G_q^{(d)} (z+1) = G_q^{(d-1)} (z) \, G_q^{(d)} (z)$,  with $G_q^{(0)} (z)=[[z]]_q$, $G_q^{(1)} (z)=\Gamma_q(z)$ and $G_q^{(2)} (z)=G_q(z)$. Here $\Gamma_q$ is the Jackson deformation  \cite{Jackson} of the Euler Gamma function. One explicit formula is:

$$
G_q(z)= \frac{(1-q)^{-\frac{1}{2} (z-2) (z-1)} \left((q;q)_{\infty
   }\right){}^{z-1} \left(\prod _{\text{s1}=0}^{\infty } \prod
   _{\text{s2}=0}^{\infty }
   \left(1-q^{\text{s1}+\text{s2}+z}\right)\right)}{\prod
   _{\text{s1}=0}^{\infty } \prod _{\text{s2}=0}^{\infty }
   \left(1-q^{\text{s1}+\text{s2}+1}\right)}
$$

When the argument is a positive integer $r$,  the quantum Barnes G-function is related as follows to the quantum superfactorial:
$
\quad G_q(r+2) =  \mathrm {Sf}_q(r)
$.
Notice that with this definition of the $G_q$, the function $ \mathrm {Sf}$, not  $\mathrm {sf}$, appears on the right hand side of the above relation.

\noindent
In the particular case $G=A_r \sim SU(r+1)$,  the results obtained for global dimensions can therefore be expressed in terms of the quantum Barnes function, since, if $r$ is a positive integer, 
$$\mathrm {sf}_q(r) = q^{-(r+1)r(r-1))/6} \, G_{q^2}(r+2) $$

In our approach, the classical and quantum Barnes functions $G(z)$ and $G_q(z)$ should be considered as associated with the Lie group $SU(r+1)$.
There should be generalizations of these functions (with complex arguments) for all simple Lie groups.

\subsubsection*{Some explicit results} 

The obtained formulae, in terms of quantum factorials of type $G$, are explicit enough, but for the convenience of the reader we give a few expressions in terms of trigonometric lines, for~$A_r$: 

\smallskip

%\paragraph{$SU(n)$}
$
\begin{array}{lll}
 \vert {\mathcal A}_k(A_1)\vert  = \frac{1}{2} (k+2) \csc ^2\left(\frac{\pi }{k+2}\right)
& {} &
%\paragraph{$SU(3)$}
 \vert {\mathcal A}_k(A_2)\vert  = \frac{3}{256} (k+3)^2 \csc ^6\left(\frac{\pi }{k+3}\right) \sec ^2\left(\frac{\pi }{k+3}\right)
\\
{} & {} & {} \\
%\paragraph{$SU(4)$}
 \vert {\mathcal A}_k(A_3)\vert  =\frac{(k+4)^3 \csc ^{12}\left(\frac{\pi }{k+4}\right) \sec ^4\left(\frac{\pi }{k+4}\right)}{16384 \left(2 \cos \left(\frac{2 \pi }{k+4}\right)+1\right)^2}
& {} &
%\paragraph{$SU(5)$}
 \vert {\mathcal A}_k(A_4)\vert  =\
\frac{5 (k+5)^4 \csc ^{20}\left(\frac{\pi }{k+5}\right) \sec ^8\left(\frac{\pi }{k+5}\right) \sec ^2\left(\frac{2 \pi }{k+5}\right)}{2^{30} \left(2 \cos
   \left(\frac{2 \pi }{k+5}\right)+1\right)^4}
 \end{array}   
$

%$$
% \vert {\mathcal A}_k(D_4)\vert  =
%\frac{(k+6)^4 \csc ^{24}\left(\frac{\pi }{k+6}\right) \sec
 %  ^8\left(\frac{\pi }{k+6}\right)}{2^{32} \left(2 \cos \left(\frac{2
  % \pi }{k+6}\right)+1\right)^6 \left(2 \cos \left(\frac{2 \pi
 %  }{k+6}\right)+\cos \left(\frac{4 \pi }{k+6}\right)+\cos \left(\frac{6
 %  \pi }{k+6}\right)+1\right)^2}
%$$
 \bigskip 
 \noindent
 For other series we only give the first few terms corresponding to levels $k=1,2,\ldots$
 
$
 \vert {\mathcal A}_k(E_6)\vert  =
 \{3,\dfrac{21}{2 \left(1-\sin \left(\frac{3 \pi
   }{14}\right)\right)},45 \left(5+2 \sqrt{5}\right),96 \left(22+15
   \sqrt{2}+4 \sqrt{58+41 \sqrt{2}}\right), \ldots\}
 $
 
 $
 \vert {\mathcal A}_k(E_7)\vert  =
\left\{2,2 \left(5+\sqrt{5}\right),21 \left(5+\sqrt{21}\right), \ldots \right\}
  $
  
   $
 \vert {\mathcal A}_k(E_8)\vert  =
  \left\{1, 4,  \epsilon, \ldots  \right\}
  $,
 $\epsilon \sim  34.64$, the highest root of {\footnotesize  $x^5-55 x^4+847 x^3-5324 x^2+14641 x-14641$}
 
  $
 \vert {\mathcal A}_k(B_2)\vert  =  \vert {\mathcal A}_k(C_2)\vert  = \{4,20, 24(2+\sqrt{3}) ,\ldots \}
  $
  
   $
 \vert {\mathcal A}_k(B_3)\vert  = \{4,28,16(4+2\sqrt{2} + \sqrt{20+14\sqrt{2}}) ,\ldots \}
  $

   $
 \vert {\mathcal A}_k(C_3)\vert  = \left\{5+ \sqrt{5}, 24( 2+\sqrt{3}) ,\ldots \right\}
  $
  
     $
 \vert {\mathcal A}_k(D_4)\vert  = \left\{ 4,32,,\ldots \right\}
  $
  
  $
 \vert {\mathcal A}_k(D_5)\vert  = \left\{ 4,40,\ldots \right\}
  $
  
  $
 \vert {\mathcal A}_k(F_4)\vert  =   \{\frac{1}{2}(5+ \sqrt{5}) , \epsilon,48(5+2 \sqrt{6}), \ldots \}
  $
  
    $
 \vert {\mathcal A}_k(G_2)\vert  = \{\frac{1}{2}(5+ \sqrt{5}) ,3 \sqrt{3 \left(5+4 \sqrt{3} \cos \left(\frac{\pi }{18}\right)+2 \cos
   \left(\frac{\pi }{9}\right)\right)}, \frac{21}{2}(5+ \sqrt{21}),\ldots \} $

\subsubsection*{Level rank duality for $A_r$}

The level-rank duality property for $A_{r} \sim SU(g=r+1)$ Lie groups was observed  in \cite{JimboMiwa}. For the coefficient $S_{00}$ of the modular matrix $S$ it implies:
$\sqrt{g} \;  S_{00}[A_{g-1},k] = \sqrt{k} \; S_{00}[A_{k-1},g] $.
In terms of global dimensions, the property reads:
$$
k \; \vert {\mathcal A}_k (A_{g-1}) \vert = g \; \vert {\mathcal A}_g (A_{k-1}) \vert
$$
It can be used to get simple enough expressions for these quantities when the rank $r$ is big and when the level $k$ is reasonably small. For instance if we want to calculate $\vert {\mathcal A}_2(A_{9})\vert$ it is easier to calculate $(10/2) \times \vert {\mathcal A}_{10}(A_1)\vert = 5 \times  24 \, (2+ \sqrt{3})$ than to use a general formula for $A_{11}$.

This property implies a duality formula for the quantum factorials $\mathrm{sf}_q(r)$ of type $A_r$,  or for the quantum Barnes function $G_q(z)$ when $z$ is a positive integer and $q$ a particular root of unity.  When $g$ and $k$ are positive integers, the two following expressions are equal:
$$
2^{\frac{1}{2} (g-1) g} \,  (g+k)^{\frac{1-g}{2}} \left( \sin \left(\frac{\pi }{g+k}\right)\right)^{\frac{1}{2} (g-1)   g}
   \mathrm{sf}_{e^{\frac{i \pi }{g+k}}} (g-1) 
$$
$$
2^{\frac{1}{2} (k-1) k} \, (g+k)^{\frac{1-k}{2}}
 \left(\sin
\left(\frac{\pi }{g+k}\right)
 \right)^{\frac{1}{2} (k-1) k}
   \mathrm{sf}_{e^{\frac{i \pi }{g+k}}} (k-1)
   $$

From an already obtained explicit expression giving for instance $\vert {\mathcal A}_k(A_1) \vert$, \ie all levels $k$ in rank $1$, we obtain immediately an explicit result giving $\vert {\mathcal A}_2(A_r) \vert = \frac{1}{4} (r+1) (r+3) \csc ^2\left(\frac{\pi }{r+3}\right)$ \ie all ranks $r$ at level $2$.
We also obtained previously an asymptotic  expression for global dimensions when the level goes to infinity. The  duality relation implies immediately, in the case $A_r$,  an asymptotic expression when the level is fixed and the rank goes to infinity:
$$
\vert {\mathcal A}_k(A_{r}) \vert  \,  \underset{{{r \rightarrow \infty}} }{\sim}\, \frac{ 1}{(2 \pi )^{(k-1) k} \,   \left( G(k+1) \right)^2} \times r^{(k^2)}
$$

%%%%%%%%%%%%%%%%%%%%%%%%

\section{Global dimensions for conformally exceptional quantum subgroups at level $k$}

\subsection{Comments} 

\subsubsection*{Purpose} 
We now suppose that the fusion category ${\mathcal A}_k(G)$ is given. We consider module-categories ${\mathcal E}$ associated with the later and we wish to calculate their global dimensions.  There are several methods to do that and everything depends on what is already known for the chosen example.
The necessary background material was summarized in section \ref{background}.
The simplest situation is when  ${\mathcal E}$ is associated with a conformal embedding: There is still a lot of work to do if we want to determine explicitly the module structure, \ie the fusion graphs (annular matrices) of  ${\mathcal E}$, the ring of its endomorphism category  ${\mathcal O}(E)=End_{\mathcal A}{\mathcal E}$ and its action, but the value of  $\vert {\mathcal E}\vert $  is particularly simple to obtain from the embedding. 
This is the simple situation that we choose to analyze in the present section.

\subsubsection*{Terminological comment.} 

We often call ``Lie group $G$ at level $k$'' the fusion category ${\mathcal A}_k(G)$. In the same way, the module-category ${\mathcal E}$ will simply be called an  ``associated quantum module''.  In what follows, we shall restrict our attention to quantum modules obtained from conformal embeddings. For a given $G$, the number of such quantum modules is finite, hence the adjective ``conformally exceptional''. 
In that case the fusion graph of ${\mathcal E}$ can be obtained from the study of the endomorphism ${\mathcal O}(E)$ and inherits a multiplication (self-fusion) from the monoidal structure of the later, for this reason ${\mathcal E}$ is called an exceptional quantum subgroup.  From the study of the same ring of endomorphisms, one sometimes  discovers one or several others quantum modules (also called exceptional) but without self-fusion. 
There may be others module-categories, not obtained from conformal embedding, but nevertheless enjoying a self-fusion compatible with the action of ${\mathcal A}_k(G)$, 
so we decide to write ``conformally exceptional quantum subgroups'' to refer to those particular quantum subgroups (existence of a compatible self-fusion) that can be constructed from conformal embeddings. This will be made more precise in a next section.

\subsection{Conformal embeddings of Lie groups at level $k$}

\subsubsection*{Definition}

The definition of conformal embeddings belongs to the lore of affine Lie algebras. However we re-write it here in a way that uses only the properties of finite dimensional Lie algebras.

{\definition{} Let $G$ be a simple Lie group, and $k$ a positive integer.  Let $J$ be a simple Lie group. 
There is a conformal embedding of $G$ at level $k$, in $J$ if the following three conditions are satisfied: 1) There is an embedding of Lie algebras $Lie(G)$ into $Lie(J)$, 2) The Dynkin index of the embedding is equal to $k$, 3) The following equality  holds: $\frac{dim(G) \times  k }{ k + g_G} = \frac{dim(J) \times 1}{ 1 + g_J}.$
Here $g_G$ and $g_J$ are the dual Coxeter numbers of $G$ and $J$. One call  $c$ the common value of the last two expressions.
If $G$ is semi-simple, not simple, we have the same definition but now $k=(k_i)$ is a multi-index, the quantity $c_i$ is defined for each simple component $G_i$ of $G$ and the equality of central charges should hold for $c=\sum c_i$.
}

\smallskip
\noindent
Remark:  In the framework of affine Lie algebras,  $c$ is interpreted as a central charge and the numbers $k$ and $1$ denote the respective levels for the affine algebras corresponding to $G$ and $J$.
The list of conformal embeddings has been known for more than twenty years, see \cite{BaisBouwknegt, KacWakimoto, SchellekensWarner}.
Here we are interested in conformally exceptional quantum subgroups of $G$ at level $k$, so there is a slight  change of perspective, in comparison with the quoted literature,
since rather than listing conformal subalgebras of $J$  we  fix $G$ and look for  ``overgroups'' $J$ such that the embedding is conformal.
Warning:  it is not difficult to find (non conformal)  embeddings $G \subset J$, and appropriate values of $k$  for which the equality of central charges is satisfied,  but where $k$ is {\sl not} the Dynkin index.  

\subsubsection*{Examples}

For the sake of illustration and in order to discuss a few examples, we shall only remind the reader what is the list of conformal embedding  for $G=SU(g)$ Lie groups.
There are three regular series and a few sporadic cases.
 
 \smallskip
 \noindent 
Regular series. The following embeddings are respectively called antisymmetric, adjoint and symmetric (take $G=SU(g)$ in all cases):
 
\smallskip
 $ 
 \begin{array}{c|c|c|c}
%G & SU(g) & SU(g) & SU(g) \\
%\hline
k &   k=g-2, &  k=g, & k=g+2, \\
{} &       g \geq 4 & g \geq 3 & g\geq 2 \\
J & SU({g(g-1)}/{2}) & Spin(g^2-1) & SU(g(g+1)/2)  \\
 \end{array}
$ 

\bigskip
 
 \noindent
Sporadic cases:

\smallskip
$
 \begin{array}{c|cc|cc|c|c|cc|c}
 G & SU(2)  & {} &  SU(3) & {}  & SU(4) & SU(6) & SU(8) & & SU(9) \\
 \hline
 k & 10 & 28 &9& 21&  8 &  6 & 1 & 10 &  1\\
J & Spin(5) & G_2  & E_6 & E_7  & Spin(20)  & Sp(10) & E_7 & Spin(70)  & E_8 \\
 \end{array}
$ 
 
 \smallskip
 \noindent
 Remarks:
 All quantum modules of $SU(2)$, $SU(3)$ and $SU(4)$ are known \cite{CIZ, DiFrancescoZuber,Ocneanu:Bariloche}.
 For $SU(2)$, their fusion graphs are classified by $ADE$ Dynkin diagrams. 
 The $E_6= {\mathcal E}_{10}(SU(2))$ and $E_8= {\mathcal E}_{28}(SU(2))$ cases respectively corresponds to the above sporadic conformal embeddings into $Spin(5)$ and $G_2$.
  The  $E_7= {\mathcal E}_{18}(SU(2))$ case (no self-fusion) can be obtained from a semi-simple conformal embeddings $SU(2)\times SU(3) \subset E_8$ followed by contraction, so it is not  ``conformally exceptional''  in our sense.
 The $D_r$ cases, $r\geq 4$ (self-fusion only when $r$ is even), can be obtained from semi-simple conformal embeddings $SU(2)\times SU(q) \subset SU(2q)$ followed by contraction, so they are not ``conformally exceptional''  either, with the exception of the smallest, $D_4= {\mathcal D}_4(SU(2))$, because it can also be constructed from the smallest symmetric regular  embedding ($g=2$). We have something similar for $g=3$ and $g=4$ because  ${\mathcal D}_3(SU(3))$ and ${\mathcal D}_2(SU(4))$
 can be  respectively constructed from the smallest adjoint and smallest antisymmetric regular embeddings. 
The regular series of adjoint conformal embeddings actually exists for all Lie groups: $G$ at level $k=g$ (its dual Coxeter number) is conformally embedded in $Spin(dim(G))$. 
 Notice finally that the embedding\footnote{Quantum symmetries of conformally exceptional quantum subgroups of $SU(4)$ are studied in \cite{RobertGil:E4E6E8ofSU4}.} 
 of $SU(4)$ in $Spin(20)$ that was flagged as sporadic from the point of view of $SU(g)$ is actually the smallest member of a regular series of embeddings for $D_r$  (indeed $A_3 \sim D_3$).
 
\subsection{From conformal embeddings to module-categories}

There are several ways to construct quantum modules or quantum subgroups from conformal embeddings. No categorical description of these constructions seems to be available so far but it is not the purpose of the present article to enter this discussion. Nevertheless we need to briefly summarize the strategy, in order to specify our own framework. Starting with some conformal embedding of $G$, at level $k$ in $J$, we consider the partition function associated with the fusion category ${\mathcal A}_1(J)$. Its matrix is the unit matrix of size $\vert {\mathcal A}_1(J)\vert \times \vert {\mathcal A}_1(J) \vert$ with lines and columns labelled by the simple objects of $J$ at level $1$. There is a known procedure (see for instance \cite{YellowBook}) to restrict the representations of $J$ to representations of $G$ existing at level $k$; this can be done, for instance in the framework of affine Lie algebras but the algorithm can be described in terms of usual Lie groups (one should be cautious since the branching rules do not usually coincide with those given by the classical theory). By restriction, the previous partition function gives a new partition function, which is block diagonal (there are $\vert {\mathcal A}_1(J)\vert $ blocks). It is now described by a block diagonal but not diagonal matrix $Z$ of size $\vert {\mathcal A}_k(G)\vert \times \vert {\mathcal A}_k(G) \vert$, with line and columns labelled by simple objects of $\vert {\mathcal A}_k(G)\vert$, and still commuting  with the action of $SL(2,\ZZ)$. The problem is then to construct a module-category ${\mathcal E}$,  with an action of a  $\vert {\mathcal A}_k(G)\vert$,  that will have $Z$ as partition function. Several techniques are available. One possibility, rather systematic but involving huge calculations, is to construct directly the endomorphism category ${\mathcal O}(E)$ from the given $Z$. This is done by using the fact that the Grothendieck ring of the later (the ``algebra of quantum symmetries'') should be a bi-module over the fusion algebra of $\vert {\mathcal A}_k(G)\vert$, and by solving the associativity equation  $(mn) x (pq) = (m(nxp)q)$ where $m,n,p,q$ refer to simple objects of $\vert {\mathcal A}_k(G)\vert$ and $x$ to simple objects of ${\mathcal O}(E)$. At the end of the day, one usually\footnote{If one finds several solutions, one should check that coherence equations for triangular cells are satisfied \cite{Ocneanu:Bariloche}.} ends up with a unique solution for  the bimodule action of the fundamental generators of the fusion ring ${\mathcal O}(E)$. This action is encoded by matrices with non-negative integer coefficients and by two graphs (because they can act from the left and from the right) for each generator;  their collection is the Ocneanu graph. 
The intersection of the connected components of the graphs associated with left and right generators,  and containing the unit of ${\mathcal O}(E)$, defines a finite subset of objects, that are called ``ambichiral vertices''.  The construction of the endomorphism algebra, that we just sketched, is fairly general and does not apply only to those examples obtained from conformal embeddings, but the later case possess specific extra features, that we list below.
 \begin{itemize}
\item The fusion graph of ${\mathcal E}$ is discovered twice by looking at the connected components (containing the unit object of ${\mathcal O}(E)$ of the left and right generators. 
\item The intersection of the two chiral copies gives defines a finite subset of objects, that are called ``modular vertices'' (as objects of ${\mathcal E}$)  in one to one correspondence with the ambichiral vertices of ${\mathcal O}(E)$. 
\item The space spanned by the vertices of the fusion graph of ${\mathcal E}$ is naturally a module over the fusion ring and over the ring of quantum symmetries, but it also inherit from the later a self-multiplication compatible with both actions.
\item From the fusion graph of ${\mathcal E}$, one can always read the annular matrices $F_n$ describing the module action of $\vert {\mathcal A}_k(G)\vert$ or, equivalently, the essential matrices $E_a$, labelled by simple objects of ${\mathcal E}$; the columns of $E_{\underline 0}$ describe induction rules. The extra feature, in the present case, is that the induction rules associated with modular objects correspond exactly to the various blocks of the partition function $Z$, \ie to the simple objects of  ${\mathcal A}_1(J)$.
 \end{itemize}
The obtained module-category ${\mathcal E}$ is what we call a ``conformally exceptional quantum subgroup'' of $G$ at level $k$.
If $G$ is semi-simple but not simple (remember that $J$ is assumed to be simple), the level $k$ is a pair of integers $(k_1,k_2)$ and there is still another construction that we want to briefly mention. Writing $G \sim G_1 \times G_2$, there is a way to perform a reduction with respect to one of the two factors (there are actually several possibilities at this point), so that one can end up with a module-category ${\mathcal E}_1$ associated with, say,  $\vert {\mathcal A}_{k_1}(G_1)\vert$.  ${\mathcal E}$ has self-fusion, but this is not always so for ${\mathcal E}_1$. Moreover the modular invariant $Z_1$ of the later is not necessarily of type 1 (not a sum of blocks).
This type of construction, \ie conformal embedding followed by reduction, is of course interesting, and leads to nice examples, but a quantum module (which may be, in some cases,  a quantum subgroup) obtained in this way, will not be called a ``conformally exceptional quantum subgroup''.
When $G=SU(g)$, rank-level duality implies the existence, for each quantum module, of a quantum module partner. 
 The rank-level dual of a conformally exceptional quantum subgroup is certainly an exceptional quantum subgroup, but it is not  always ``conformally exceptional''  in our sense, because its construction may require a semi-simple, non simple, conformal embedding followed by some kind of non-diagonal reduction.
 
 \subsection{On the global dimension of exceptional quantum subgroups}

 \subsubsection*{Constraints from induction-restriction rules}
 
It is convenient to think of $\mathcal{A}_k(G)/\mathcal{E}$ as a homogenous space, both discrete and quantum.
Classically, the dimension of the space $\Gamma_a$ of vector valued functions defined on the quotient of a finite group $A$ by a subgroup $E$, and valued in a vector space $a$, can be calculated either trivially as $dim(a)\times \vert A/E \vert$ or non trivially, when $a$ is a representation space for $E$,  by decomposing this space of functions  into a sum of irreducible representations of $A$; this applies in particular to the case where $a$ is the trivial representation ${\underline 0}$ of $E$, so that we obtain a decomposition of the space ${\mathcal F}$ (actually the algebra) of complex valued functions over $A/E$, and $\vert {\mathcal F} \vert = dim(\Gamma_{\underline 0})$. Here ${\mathcal F}$ refers both to a discrete set and to its algebra of functions.
In our situation, there is an induction-restriction functor \cite{Ostrik}, so that the above still makes sense. 
In this analogy, vertices of the fusion graphs of $\mathcal{E}$ do not only label simple objects $a$ of the later  but also spaces $\Gamma_a$ of sections of quantum vector bundles  which can be decomposed, using induction, into simple objects $n$ of $\mathcal{A}_k(G)$: We write $\Gamma_a = \sum _{n  \uparrow \Gamma_a} \,  n$.
Like in the classical situation, the (quantum) dimension of $\Gamma_a$ can be calculated in two ways :
$$qdim(\Gamma_a)=\sum _{n  \uparrow \Gamma_a} \,  qdim(n) \quad  \text{and} \quad  qdim(\Gamma_a)=qdim(a) \times \vert {\mathcal A}_k(G)/{\mathcal E}\vert$$
There is a special object of ${\mathcal E}$, denoted ${\underline 0}$, with quantum dimension is equal to $1$, and whose space of sections ${\mathcal F} \doteq\Gamma_{\0}$ plays the role of an algebra of functions over a ``quantum  space'', it is actually a Frobenius algebra.
When $a={\0}$, the previous equation gives:  $$\vert {\mathcal F} \vert = \sum _{n  \uparrow \Gamma_{\0}} \,  qdim(n) =  \vert {\mathcal A}_k(G)/{\mathcal E}  \vert $$
If ${\mathcal E}$ is associated with a conformal embedding, ${\underline 0}$ plays the role of unity for self-fusion. As already mentioned the modular invariant is then a sum of blocks (we are in the so-called ``type I'' case) and  the decomposition of ${\mathcal F}$ into simple objects can be read from the first modular block of the partition function $Z$, or from the first line of the associated matrix. 

\subsubsection*{Constraints on quantum dimensions from the modular invariant}

Call $Z$ the partition function of ${\mathcal E}$.
We shall need later the known identity:
$$\sum_{m,n}   qdim(m) \,  Z_{m,n} \, qdim(n) = \sum_n  qdim(n)^2 \,   = \vert {\mathcal A}_k(G) \vert $$
Indeed, from the expression of quantum dimensions of simple objects $m,n$ of ${\mathcal A}_k(G)$ in terms of matrix elements of the modular generator $S$, 
the left hand side reads $\sum_{m,n} S_{0m} Z_{mn} S_{0n}/S_{00}^2$. But $S_{0m}$ is real and $S$ is symmetric, so that the l.h.s. reads $(S Z S^\dag)_{00}/S_{00}^2$, or  $(Z S S^\dag)_{00}/S_{00}^2$
since $SZ=ZS$ from modular invariance of $Z$. Unitarity of $S$ together with the normalization condition $Z_{00}=1$ imply that the l.h.s. is just $1/S_{00}^2$, hence the result.
Incidentally, this relation can be used to provide a check of the correctness of the calculated modular invariant $Z$.
 
 \smallskip
 The simple objects of  $ {\mathcal A}_k(G)$ labelling non-zero diagonal entries  of the modular invariant matrix $Z$ define the family (possible multiplicity) ${\mathcal E}xp$ of generalized exponents of ${\mathcal E}$. When $G=SU(2)$, one recovers the usual exponents for $ADE$ Lie groups.
 In those cases obtained from conformal embeddings,  $Z$ is a sum of blocks (type $I$), we have therefore a partition of the set of exponents. The different parts of this partition are in one to one correspondence with the modular vertices of ${\mathcal E}$.
 
\subsubsection*{A simple formula for $\vert {\mathcal E} \vert$}
As usual, $m,n,\ldots$ label simple objects of $ {\mathcal A}$ (we suppress below the reference to the group $G$ and to the level $k$), and $a,b,\ldots$, simple objects of  ${\mathcal E}$.
The formula that we shall obtain for $\vert{\mathcal E}\vert$ relies on the following steps:
\begin{enumerate}
\item By definition, $\vert {\mathcal A} \vert = \sum_n \, qdim(n)^2$,  $\vert {\mathcal E} \vert = \sum_a \, qdim(a)^2$, and $\vert {\mathcal F} \vert = \vert {\mathcal A} \vert /\vert {\mathcal E} \vert$.
\item For the special object $0$ of  $ {\mathcal E}$, we have $ {\mathcal F} = \oplus_{n  \uparrow \Gamma_0}  \, n$. Then  $\vert {\mathcal F} \vert = qdim(\Gamma_0) =  \sum_{n  \uparrow \Gamma_0}  \, qdim(n)$.
\item For simple objects of  $ {\mathcal E}$ we have, more generally, $qdim(a) \times \vert {\mathcal F} \vert  =  qdim(\Gamma_a) =  \sum_{n  \uparrow \Gamma_a}  \, qdim(n)$.
\item From unitarity of the modular matrix $S$ (see above),  $\sum_{m,n}   qdim(m) \,  Z_{m,n} \, qdim(n) = \sum_n  qdim(n)^2 \,   = \vert {\mathcal A} \vert $, where $Z=(Z_{m,n})$ is the matrix defining the partition function.
\item For conformally exceptional quantum subgroups, $Z$ is a sum of blocks (type $I$) labelled by modular vertices $a \in J$. 
Moreover, for $a \in J$ and defining $Z_{m,n}^{(a)}$ by $(qdim(\Gamma_a))^2 = \sum_{m,n  \uparrow \Gamma_a} qdim(m) \, Z_{m,n}^{(a)} \, qdim(n)$, the l.h.s. of the relation obtained in step 4 can be re-written $\sum_{a\in J}  \sum_{m,n  \uparrow \Gamma_a}  qdim(m) \, Z_{m,n}^{(a)} \, qdim(n)$. Therefore, $\sum_{a \in J} \, (qdim(\Gamma_a))^2  = \vert {\mathcal A} \vert$.
\item From the relation obtained in step 3, and defining $\vert {\mathcal J} \vert =  \sum_{a \in J} \, qdim(a)^2$, the equality obtained in step 5 reads $\vert {\mathcal J} \vert \times \vert {\mathcal F} \vert^2 = \vert {\mathcal A} \vert$. But $\vert {\mathcal F} \vert = \vert {\mathcal A} \vert /\vert {\mathcal E} \vert$,  therefore $\vert {\mathcal A} \vert = \vert {\mathcal E} \vert^2/   \vert {\mathcal J} \vert$.
\end{enumerate}
This last equality has actually a range of validity wider than the one discussed in this article and relies on the fact  \cite{Ocneanu:talks} that the ring of  ${\mathcal O}({\mathcal E}) = End_{{\mathcal A}_k(G)}({\mathcal E)}$ is isomorphic to  ${\mathcal E}_L \otimes {\mathcal E}_R / {\mathcal J}$, \ie to the tensor product of left and right chiral subalgebras over the ambichiral subalgebra.
The previous equality relating global dimensions then follows from $\vert{\mathcal A}_k(G)\vert = \vert {\mathcal O}({\mathcal E})\vert$ and from the fact that, in the case of a quantum module ${\mathcal E}$  measuring a conformal embedding, the following properties hold:
1) Both ${\mathcal E}_L$ and ${\mathcal E}_R$ are isomorphic with ${\mathcal E}$ (which has indeed self-fusion : it is quantum subgroup), 2) The set $J$ of  modular vertices of ${\mathcal E}$ can be identified with the simple objects of the subcategory  ${\mathcal J}$ of  ${\mathcal O}({\mathcal E})$, 3) ${\mathcal J}$ itself, with global dimension $\vert {\mathcal J} \vert$ (called ambichiral dimension of $\vert {\mathcal E} \vert$) can be identified with the fusion category, at level $1$, of the Lie group in which $G$ is conformally embedded. This explains the notation $J$, that we chose to denote it, since ${\mathcal J} = {\mathcal A}_1(J)$.  Details about the above more general result describing the structure of the ring of  ${\mathcal O}({\mathcal E})$, in the general case,  are unfortunately unavailable, so we decided to present a simpler argument since we are anyway only concerned in this paper with calculating global dimensions for conformally exceptional quantum subgroups.
Notice that $ \sqrt{ \vert {\mathcal A}_k(G)}$ coincides with the inverse $1/S_{00}$ of the first matrix element of the modular $S$ matrix.  The conclusion is that  $\vert {\mathcal E} \vert $ can be obtained without having to calculate independently the quantum dimensions of the simple objects (even without knowing how many there are). Let us summarize: 

 \label{formulaforE}
 
 {\theorem {} Consider a conformal embedding of the simple or semi-simple Lie group $G$, at level $k$, in the simple Lie group $J$.  This embedding is associated (or ``measured by'') a module-category 
  ${\mathcal E}$, with an action of the fusion category ${\mathcal A}_k(G)$. Its global dimension is given by the formula:
  $$  \vert {\mathcal E} \vert = \sqrt {\vert {\mathcal A}_k(G)\vert \times   \vert {\mathcal J} \vert } $$
  where   $\vert {\mathcal J} \vert $ is the global dimension of  ${\mathcal A}_1(J)$.
} 
 
 \subsubsection*{Discussion (illustrated with an example of type $SU(2)$)}
 
 Calculating the global dimension $\vert {\mathcal E} \vert$ can be done in many ways, the best strategy depends about what is actually known on the example at hand. 
 We shall mention  several possibilities and  illustrate the discussion with the example of $E_8= {\mathcal E}_{28}(SU(2))$ associated with the fusion category $A_{29} = {\mathcal A}_{28}(SU(2))$. The later has $29$ simple objects $n$, with quantum dimensions $qdim(n) = [n]_q$,  $q = \exp(2 i \pi/30)$, for $n=1,2,\ldots, 29$. The Perron-Frobenius norm is  $[2]_q = 2 \cos \frac{\pi}{30}$.
 The global dimension  $\vert A_{29}\vert = 30(12+5\sqrt{5} + \sqrt{3(85+38\sqrt{5})} )$ can be obtained by summing squares, or from general formulae (first section).
\begin{itemize}
\item Possibility 1: The fusion graphs are known (or equivalently, the adjacency matrices). In our example, the adjacency matrix is $F_1=C-2 \one$ where $C$ is the Cartan matrix of $E_8$. The Perron-Frobenius eigenvalue is $[2]_q$. The quantum dimensions $qdim(a)$ of the $8$ simple objects are given by the $F_1$ eigenvector  $([1]_q, [2]_q, [3]_q, [4]_q, [5]_q, [7]_q/[2]_q, [5]_q/[3]_q; [5]_q/[2]_q)$. The global dimension is obtained from its definition : $\vert E_8 \vert = \sum_a \, qdim(a)^2 = \dfrac{1}{2}(15(3+\sqrt{5}) + \sqrt{30(65+29\sqrt{5})})$.
 
\item  Possibility 2:  Quantum dimensions can also be read from the induction rules (given by the columns of the essential matrix $E_{\underline 0}$).
One evaluates, in turns,  $qdim(\Gamma_a)=\sum _{n  \uparrow \Gamma_a} \,   qdim(n)$, in particular $\vert {\mathcal F} \vert = qdim(\Gamma_{\underline 0})$, then $qdim(a) = qdim(\Gamma_a)/qdim(\Gamma_{\underline 0})$. Finally $\vert {\mathcal E} \vert$ is obtained by summing squares. Actually, if we are only interested in the later result, it is enough to calculate $\vert {\mathcal F} \vert$ since
$\vert {\mathcal E} \vert =  \vert A_{k}\vert / \vert {\mathcal F} \vert$. In our case, $\vert {\mathcal F} \vert = [1]_q+[11]_q+[19]_q+[29]_q = \dfrac{1}{2}(3(5+\sqrt{ 5})+\sqrt{150 + 66 \sqrt{5}})$. One recovers
 $\vert E_8 \vert = \vert A_{29}\vert / \vert {\mathcal F} \vert$. 
 
 \item If the partition function $Z$  is known and if it  is block diagonal (like in our example), it determines the induction rules for the set $J$ of modular points, in particular for the origin. This is enough to determine $\vert {\mathcal E} \vert$ since only the first block is needed. Calculating the quantum dimensions of the modular points, via the formula $qdim(a) = qdim(\Gamma_a)/qdim(\Gamma_{\underline 0})$ provides a check since $\vert {\mathcal E} \vert$ can also be obtained from $ \vert {\mathcal A}_k \vert$ and $ \vert {\mathcal J} \vert = \sum_{a \in J} \, qdim(a)^2$. In our case, $Z = (\chi_1 + \chi_{11} + \chi_{19} + \chi_{29})^2 + (\chi_{7} + \chi_{13} + \chi_{17}+\chi_{23})^2$. Here, the labels $n$ denoting the  simple objects of $A_{29}$ have been shifted by $1$. 
 There are two modular points (the two extremal vertices of the longest branches of $E_8$) with q-dimensions $[1]$ and $[5]_q/[3]_q$, so that $ \vert {\mathcal J} \vert = \dfrac{1}{2}(5+\sqrt{5})$. 
 One recovers   $\vert E_8 \vert = \sqrt{ \vert A_{29}  \vert \vert {\mathcal J} \vert }$. 

\item If we a priori know that the chosen example is defined by a conformal embedding, like in this example which can indeed be obtained as the conformally exceptional quantum subgroup measuring the embedding of $SU(2)$ at level $28$ in $G_2$, the easiest is to determine $\vert {\mathcal F} \vert$ from the later ``overgroup''.  Here it is equal to  the global dimension $\vert {\mathcal A}_1(G_2) \vert = (5+ \sqrt{5})/2$. One again recovers the value  $\vert E_8 \vert$.

\end{itemize} 
The reader willing to play with this example may use the induction rules given below:  they are non-zero entries of the essential matrix $E_{\underline 0}$. Lines are labelled by simple objects of $A_{29}$ and columns by simple objects of $E_8$ (the vertex belonging to the shortest branch is at the end. Modular vertices are in position $1$ and $7$. For typesetting reasons we display it horizontally (so we give the transpose of  $E_{\underline 0}$).

\smallskip
\noindent
{\tiny
$
\left(
\begin{array}{ccccccccccccccccccccccccccccc}
 1 & . & . & . & . & . & . & . & . & . & 1 & . & . & . & . & . & . & . & 1 & . & . & . & . & . & . & . & . & . & 1 \\
 . & 1 & . & . & . & . & . & . & . & 1 & . & 1 & . & . & . & . & . & 1 & . & 1 & . & . & . & . & . & . & . & 1 & . \\
 . & . & 1 & . & . & . & . & . & 1 & . & 1 & . & 1 & . & . & . & 1 & . & 1 & . & 1 & . & . & . & . & . & 1 & . & . \\
 . & . & . & 1 & . & . & . & 1 & . & 1 & . & 1 & . & 1 & . & 1 & . & 1 & . & 1 & . & 1 & . & . & . & 1 & . & . & . \\
 . & . & . & . & 1 & . & 1 & . & 1 & . & 1 & . & 1 & . & 2 & . & 1 & . & 1 & . & 1 & . & 1 & . & 1 & . & . & . & . \\
 . & . & . & . & . & 1 & . & 1 & . & . & . & 1 & . & 1 & . & 1 & . & 1 & . & . & . & 1 & . & 1 & . & . & . & . & . \\
 . & . & . & . & . & . & 1 & . & . & . & . & . & 1 & . & . & . & 1 & . & . & . & . & . & 1 & . & . & . & . & . & . \\
 . & . & . & . & . & 1 & . & . & . & 1 & . & . & . & 1 & . & 1 & . & . & . & 1 & . & . & . & 1 & . & . & . & . & .
\end{array}
\right)
$
}

 \subsection{Examples}
 
  \subsubsection*{Global dimensions for Lie groups at level $1$}
 
 In order to make use of the formula  \ref{formulaforE} for $\vert {\mathcal E} \vert$, we need to know, for each Lie group, the global dimensions $\vert {\mathcal A}_1(G) \vert$, that we just call $\vert {\mathcal A}_1 \vert$ here. We give below the list of representations existing when $k=1$,  the corresponding vector $Q$ of  quantum dimensions, and the global dimension $\vert {\mathcal A}_1 \vert$. It can be calculated as $\vert Q \vert^2$ or from the general results of the first part.
As usual,  we call basic representations those fundamental representations of smallest classical dimension:  for $E_6$ they are of highest weight $(1,0,0,0,0;0)$ and $(0,0,0,0,1;0)$, both of classical dimension $26$; for $G_2$, this is the representation of highest weight $(0,1)$ and classical dimension $7$; etc.
\begin{itemize}
\item Type $A_r$ :  The trivial and all the fundamental. $Q=\{1,1,\ldots 1\}$.  $\vert {\mathcal A}_1 \vert= r+1$.
\item Type $D_r$ : The trivial, the two half-spinorial and the vectorial. $Q=\{1,1,1,1\}$. $\vert {\mathcal A}_1 \vert= 4$.
\item Type $E_6$ : The trivial and the two basic representations. $Q=\{1,1,1\}$. $\vert {\mathcal A}_1 \vert= 3$.
\item Type $E_7$ : The trivial and the  basic representation. $Q=\{1,1\}$. $\vert {\mathcal A}_1 \vert= 2$.
\item Type $E_8$ : Only the trivial representation. $Q=\{1\}$. $\vert {\mathcal A}_1 \vert= 1$.
\item Type $B_r$ : The trivial, the spinorial and the vectorial.  $Q=\{1,\sqrt 2,1\}$. $\vert {\mathcal A}_1 \vert= 4$.
\item Type $F_4$ : $Q=\{1, \frac{1}{2}(1 + \sqrt 5)\}$. $\vert {\mathcal A}_1 \vert=  \frac{1}{2}(5 + \sqrt 5)$.
\item Type $G_2$ : The trivial and the  basic representation. $Q=\{1, \frac{1}{2}(1 + \sqrt 5)\}$. $\vert {\mathcal A}_1 \vert=  \frac{1}{2}(5 + \sqrt 5)$.
\item Type $C_r$ :   The trivial and all the fundamental. No simple closed formulae in general. We give results for low rank.  $C_2 $ is  like $ B_2$.  Type $C_3$:  $Q=\{1,1,(1+\sqrt 5)/2,(1+\sqrt 5)/2\}$.   $\vert {\mathcal A}_1 \vert=  (5 + \sqrt 5)$.
Type $C_4$:  $Q=\{1,1, \sqrt 3, 2 , \sqrt 3 \}$.   $\vert {\mathcal A}_1 \vert= 12$.
\end{itemize}

Using both the above results, giving $\vert {\mathcal A}_1(J) \vert$,  and the  general formulae of the first part, expressing global dimensions of $\vert {\mathcal A}_k(G) \vert$ as quantum superfactorials of type $G$, we can now find explicit results for global dimensions of $\vert {\mathcal E} \vert$ by using  formula \ref{formulaforE}. 

\subsubsection*{$ADE$ trigonometric identities for quantum subgroups of type $SU(2)$}

For conformally exceptional quantum subgroups of type $SU(2)$, we recover the known results:

$
\begin{array}{lccc}
D_4: &  \quad \vert {\mathcal A}_4 \vert = 12 \quad \text {and}   \quad \vert {\mathcal J} \vert = 3 , \quad  \text {so} \quad & \vert D_4 \vert =  6, \\ 
E_6: & \quad \vert {\mathcal A}_{10} \vert = 24(2+\sqrt{3}) \quad \text {and}  \quad \vert {\mathcal J} \vert = 4 , \quad  \text {so} \quad & \vert E_6 \vert =   4(3+\sqrt{3}), \\
E_8: & \quad \vert {\mathcal A}_{28} \vert = 30(12+5\sqrt{5} + \sqrt{3(85+38\sqrt{5})} ), &   \\
{} &  \text {and}  \, \vert {\mathcal J} \vert = (5+\sqrt{5})/2, \quad  \text{so} \quad & \vert E_8 \vert   =\dfrac{1}{2}(15(3+\sqrt{5}) + \sqrt{30(65+29\sqrt{5})}).
\end{array}
$

Given ${\mathcal E}$, and its partition function $Z$, it should be clear that any identity of the type $Q.Z.Q=\vert {\mathcal A}_k(G) \vert$, where $Q$ is the dimension vector of ${\mathcal A}_k(G)$,
can be interpreted as a trigonometric identity.
For $SU(2)$ we have a family of $ADE$ trigonometric identities since its module-categories are classified by simply laced Dynkin diagrams. 
For any quantum module of $SU(2)$ at level $k$,  setting $\kappa=k+2$, and after simplification by a common denominator, those identities read:
 $$  \sum_{m,n=1}^{r}   (Z)_{m,n} \, \sin(m\, \pi/\kappa) \,    \sin(n \, \pi/\kappa)  =  \kappa/2$$
At level $28$ for example, \ie at the altitude $\kappa = 28+2=30$ interpreted here as a dual Coxeter number, we have three identities, one for the fusion category itself, $A_{29} = {\mathcal A}_{28}(SU(2))$, with $Z=1$, one for the quantum module $D_{16} =  {\mathcal D}_{28}(SU(2))$, and one for the (conformally) exceptional quantum subgroup $E_8 =  {\mathcal E}_{28}(SU(2))$. 
For instance one obtains in the case $E_8$  (notice the appearance of exponents of $E_8$,  and the two modular blocks reminiscent of the conformal embedding of $SU(2)$ at level $28$ in $G_2$, since ${\mathcal A}_1(G_2)$ has only two simple objects):
{\scriptsize
\begin{eqnarray*}
 \Big\{{\sin(7 \frac{\pi}{30})} + {\sin(13 \frac{\pi}{30})} + {\sin(17 \frac{\pi}{30})}  + {\sin(23 \frac{\pi}{30})}\Big\}^2 &+& \Big\{{\sin(1 \frac{\pi}{30})} + {\sin(11\frac{\pi}{30})} + {\sin(19 \frac{\pi}{30})}  + {\sin(29 \frac{\pi}{30})}\Big\}^2  = 15
\end{eqnarray*}
  }

\subsubsection*{Regular conformal embeddings for unitary groups}
To illustrate the formula expressing $\vert {\mathcal E} \vert $ in terms of $\vert{\mathcal A}_k(G) \vert $ and $\vert {\mathcal J} \vert $  , we give the global dimensions of the first few  exceptional quantum subgroups measuring  regular conformal embeddings of type $SU(g)$ : 
\begin{itemize}
\item Regular antisymmetric series ($k=g-2, g=4,5,\ldots$). Using $\vert{\mathcal J}\vert= g(g-1)/2 $, one finds  $\vert {\mathcal E}_{g-2}(SU(g)) \vert =
\left\{12, 20 \left(2+\sqrt{2}\right),60 \left(5+2 \sqrt{5}\right),504 \left(7+4 \sqrt{3}\right), \ldots\right\}$
\item Regular adjoint series  ($k=g, g=3,4,5,\ldots$). Using $\vert{\mathcal J}\vert= 4$, one finds  $\vert {\mathcal E}_{g}(SU(g)) \vert = 
\left\{12,16 \left(2+\sqrt{2}\right),40 \left(5+2 \sqrt{5}\right),288 \left(7+4 \sqrt{3}\right), \ldots \right\}$
\item Regular symmetric series  ($k=g+2, g=2,3,4,5,\ldots$). Using $\vert{\mathcal J}\vert= g(g+1)/2$, one finds  $\vert {\mathcal E}_{g+2}(SU(g)) \vert = 
\left\{6,12 \left(2+\sqrt{2}\right),40 \left(5+2 \sqrt{5}\right),360 \left(7+4 \sqrt{3}\right), \ldots \right\}$
\end{itemize}
 
Members of the regular adjoint series are rank-level self-dual, but the rank-level dual of a conformally exceptional quantum subgroup of $SU(g-2)$ belonging to the regular symmetric series  is a  conformally exceptional quantum subgroup of $SU(g)$ 
belonging to the regular antisymmetric series, and reciprocally.  The ratio of their dimensions is: 
 $$\frac{\vert {\mathcal E}_{g-2}(SU(g)) \vert}{\vert {\mathcal E}_{g}(SU(g-2 )) \vert}  = \frac{g}{g-2}.$$
Indeed this ratio is obtained as the product  of two terms 
$\sqrt{  \vert {\mathcal A}_{g-2}(SU(g)) \vert /  \vert {\mathcal A}_{g}(SU(g-2 )) \vert }$ and $ \sqrt {\frac {g(g-1)/2}{(g-2)(g-1)/2} }$, both giving a contribution $\sqrt \frac {g}{g-2}$, hence the result.

 \subsubsection*{The biggest exceptional of the biggest exceptional }
 Let us conclude this section by giving the global dimension of the biggest conformally exceptional subgroup of the exceptional Lie group $E_8$. It occurs at level $30$ (an altitude $\kappa = 30+30=60$) and it measures the adjoint embedding of $E_8$ in $D_{124}\sim Spin(248)$.  The associated fusion-category ${\mathcal A}_{30}(E_8)$ has $20956$ simple objects and we did not calculate the number of simple objects of its
module-category ${\mathcal E}_{30}(E_8)$. Nevertheless, with $q = \exp \frac{i \pi}{60}$ and $\vert {\mathcal A}_1(D_{124}) \vert = 4$, one finds immediately the global dimension :  
 $$
\vert {\mathcal E}_{30}(E_8)\vert =  2 \times \frac{60^4}{2^{120} \;  [1]!_q [7]!_q [11]!_q [13]!_q [17]!_q [19]!_q [23]!_q [29]!_q \;  {\sin}\left[\frac{\pi }{60}\right]^{120}} \sim 5.57902\times 10^{22}
 $$

\section{Geometrical and physical considerations: Chern-Simon theory and strings}

As it is known since \cite{WittenCS95},  a Chern-Simons gauge theory in three dimensions can be viewed as a string theory.
The Chern - Simon partition function, with simply connected gauge group $G$, on the closed orientable $3$-manifold $M$, with $k$ an integer, 
is formally given by the following functional integral over the space of connections defined on a $G$ principal bundle over $M$:
 $$Z_{CS}[M,G,k] = \int {\mathcal D} A \; \exp[  i k \times  \frac{1}{4 \pi} \int_M d^3x \, Tr(A dA + \frac{2}{3} A^3) ]$$

Using a normalization $Z_{CS}[S^2 \times S^1,G,k] =1$, it was shown  in \cite{WittenCS} that, with $M = S^3$, $Z_{CS}[S^3,G,k] = S_{00}$.
where $S$ is the matrix representing the modular generator $\tau \mapsto -1/\tau$ of $SL(2,\ZZ)$ for $G$ at level $k$.
Still in \cite{WittenCS}, it was shown that $Z[T^3, SU(N), k]$ counts the number of integrable irreducible highest-weight representations of the affine Kac-Moody algebra at level $k$, \ie with another terminology, the number of simple objects in the fusion category ${\mathcal A}_k(SU(N))$.  This last number, equal to $\frac{1}{k \, B(k,N)}$, where $B$ is the Euler Beta function, was calculated by \cite{Periwal}.
The value of  $Z[T^3, SU(N), k]$ was calculated in \cite{KacWakimoto} : 
$$Z_{CS}[S^3,SU(N),k]  =   (N+k)^{-N/2} \,  \sqrt{\frac{N+k}{N}} \,  \prod _{j=1}^{N-1} 2^{N-j} \sin ^{N-j}\left(\frac{\pi  j}{N+k}\right)$$
Of course, the result can be expressed in terms of  global dimensions  for fusion categories of $SU(N)$ at level $k$ since  $\vert {\mathcal A}_k(SU(N)) \vert =1/S_{00}^2$.
The quantum Lie superfactorials of type $G$ defined in the first part of this paper, for all simple Lie groups $G$ can, in turn, be used to evaluate explicitly Chern-Simon partition functions on $S^3$.

As already mentioned in the introduction, simple objects of module-categories  ${\mathcal E}$  associated with ${\mathcal A}_k(G)$ can be interpreted in terms of boundaries for conformal field theories (Wess-Zumino-Witten models), simple objects of ${\mathcal O}({\mathcal E})=End_{{\mathcal A}_k(G)}{\mathcal E}$ in terms of defects, and presumably objects of $Hom_{{\mathcal A}_k(G)}({\mathcal E}_1,{\mathcal E}_2)$ in terms of interfaces, but when ${\mathcal E} \neq {\mathcal A}_k(G)$, we are not aware of any precise interpretation, in terms of differential geometry, or in terms of string theory, for the global dimensions.

\section*{Appendix}
\subsection*{Tables of long indices and Coxeter numbers}
For convenience we remind the reader the  values of Coxeter numbers $\gamma$, dual Coxeter numbers $g$, and  the ``long index''  for all Lie groups. The quantity $\Delta$, used in the text, is the determinant of the fundamental quadratic form,  \ie the inverse of the long index.

\smallskip
\begin{center}
$
\begin{array}{cccccccccc}
{} & A_r & B_r& C_r& D_r & E_6 & E_7& E_8 & F_4& G_2 \\
\gamma : & r+1 &2r & 2r &2r-2 & 12& 18 & 30 & 12& 6 \\
g :  &  r+1&2 r -1 &r+1 &2r-2 &12 &18 &30 & 9& 4 \\
\Delta^{-1}: & r+1 & 4 & 2^r & 4 & 3 & 2 & 1 & 4 & 3
\end{array}
$
\end{center}

% \vfill \eject

  \end{document}